\documentclass[12pt]{article}
\usepackage{epsfig,amsmath,amssymb,stmaryrd,enumerate}
\usepackage{graphicx}

\newtheorem{theorem}{Theorem}
\newtheorem{lemma}{Lemma}

\advance \textheight by \topskip
\advance \textheight by 1pt
\makeatother

\begin{document}

\title {\hbox{\normalsize }\hbox{}
On Rank Two Algebro--Geometric Solutions of an Integrable Chain \thanks{The work was supported by RSF (grant 14-11-00441)}}

\author {Gulnara S.~Mauleshova, \ Andrey E.~Mironov}

\date{}
\maketitle

\begin{abstract}
In this paper we consider a differential--difference system which is equivalent to the commutativity condition of two differential--difference operators. We study the rank two algebro--geometric solutions of this system.
\end{abstract}

\maketitle

\section{Introduction and Main Results}


I.M.~Krichever and S.P.~Novikov~\cite{KN2,KN1} proved the existence of the rank $l>0$ algebro--geometric solution of the Kadomtsev--Petviashvili (KP) equation and the Toda chain. For such solutions the common eigenfunctions of auxiliary commuting operators (differential in the case of KP or difference in the case of the Toda chain) form a rank $l$ vector bundle over the affine part of the spectral curve. They also proved that in the case of the rank two solutions of KP corresponding to elliptic spectral curves there is a remarkable separation of variables. Such solutions are expressed through the solutions of the 1+1 Krichever--Novikov (KN) equation (see (\ref{eq6}) below) and the solutions of an ODE~\cite{KN2} (see also~\cite[ formula (22)]{Mir}).

In this paper we study the rank two algebro--geometric solutions of the following equation
\begin{equation}
\label{eq1}
[\partial_y - T - f_n(x,y), \partial_x - b_n(x,y) T^{-1}- d_n(x,y) T^{-2}] = 0,
\end{equation}
where $f_n, b_n, d_n$ are the 4--periodic functions, $f_{n+4} = f_n, b_{n+4} = b_n$,\linebreak $d_{n+4} = d_n.$

Equation (\ref{eq1}) is equivalent to the following 4--periodic chain
\begin{gather}
\label{eq2}
f_{n,x}(x,y) - b_n(x,y) + b_{n+1}(x,y) = 0,\\
\label{eq3}
f_{n-2}(x,y) - f_n(x,y) + \frac{d_{n,y}(x,y)}{d_{n}(x,y)} = 0,\\
\label{eq4}
f_{n-1}(x,y) - f_n(x,y) + \frac{b_{n,y}(x,y)}{b_{n}(x,y)} +\frac{d_n(x,y)-d_{n+1}(x,y)}{b_n(x,y)} = 0.
\end{gather}
Recall that if two difference operators
$$
L_k=\sum^{K_+}_{j=-K_-} u_j(n)T^j, \quad
L_m=\sum^{M_+}_{j=-M_-}v_j(n)T^j,  \quad n \in \mathbb{Z},
$$
where $T$ is the shift operator, $T \psi_n = \psi_{n+1},$ commute, then there is a polynomial $R(z,w)$ such that $R(L_k,L_m)=0.$ The {\it spectral curve} $\Gamma$ is defined by the equation $R=0.$ The spectral curve parametrizes  the common eigenvalues, i.e. if $\psi_n$ is a common eigenfunction of $L_k, L_m$
$$
L_k \psi_n = z \psi_n, \qquad L_m \psi_n = w \psi_n,
$$
then $P=(z,w) \in \Gamma.$ The {\it rank $l$} of the pair $L_k, L_m$ is
$$
l={\rm dim}\{\psi_n:L_k\psi_n=z\psi_n, \quad L_m\psi_n=w\psi_n\}
$$
for the general $P=(z,w) \in \Gamma.$ The maximal commutative ring of difference operators is isomorphic to the ring of meromorphic functions on a spectral curve (a closed Riemann surface) with poles
$q_1,\ldots,q_s.$ Such operators are called {\it $s$--point operators}.  In the case of the rank one operators, the eigenfunctions (Baker--Akhiezer functions) can be found explicitly in terms of theta--functions of the Jacobi variety of spectral curves, and coefficients of such operators can be found using eigenfunctions. The case of higher rank is very complicated (higher rank Baker--Akhiezer functions are not found). The one--point rank two operators in the case of elleptic spectral curves were found in~\cite{KN2}. The one--point rank two operators in the case of the hyperelliptic spectral curve
\begin{equation}
\label{eq5}
w^2=F_g(z)=z^{2g+1}+c_{2g}z^{2g}+c_{2g-1}z^{2g-1}+\cdots+c_0
\end{equation}
were studied in~\cite{MM1}. In particular, examples of such operators were found for an arbitrary $g>1:$
\begin{enumerate}
\item[1)] the operator
$$
L^{^\sharp}_4=(T+(r_3n^3+r_2n^2+r_1n+r_0)T^{-1})^2+g(g+1)r_3n, \qquad r_3\ne 0
$$
commutes with a difference operator $L^{^\sharp}_{4g+2},$

\item[2)] the operator
$$
L^{\flat}_4=(T+(r_1\cos(n)+r_0)T^{-1})^2-4r_1\sin(\frac{g}{2})\sin(\frac{g+1}{2})\cos(n+\frac{1}{2}),\qquad r_1\ne 0
$$
commutes with a difference operator $L^{\flat}_{4g+2}.$
\end{enumerate}

Following~\cite{KN2,KN1}, we call the solution $f_n, b_n, d_n$ of (\ref{eq2})--(\ref{eq4}) the {\it algebro--geometric solution of rank two,} if there are one--point rank two commuting difference operators
$$
L_4=\sum^{2}_{i=-2}u_j(n,x,y)T^j, \qquad L_{4g+2}=\sum^{2g+1}_{i=-(2g+1)}v_j(n,x,y)T^j
$$
commuting with
$\partial_x - b_n(x,y) T^{-1}- d_n(x,y) T^{-2}$ and $\partial_y - T - f_n(x,y).$

In the next theorem we show that in the case of an elliptic spectral curve given by the equation
\begin{equation}
\label{eq6}
w^2 = F_1(z) = z^3 +c_2 z^2+c_1 z+c_0
\end{equation}
there is a separation of variables for rank two genus one solutions (similar to KP) of (\ref{eq2})--(\ref{eq4}).
\begin{theorem}
Let $\gamma_n=\gamma_{n}(x)$ and $\wp(y)$ satisfy the equations
\begin{gather}
\label{eq7}
\gamma'_n = \frac{F_1(\gamma_n)(\gamma_{n-1}-\gamma_{n+1})}{(\gamma_{n-1}-\gamma_n)(\gamma_{n}-\gamma_{n+1})},\\
\label{eq8}
(\wp'(y))^2 = F_1(\wp(y)),
\end{gather}
and $\gamma_{n+4}=\gamma_n,$ then
\begin{gather*}
b_n(x,y) = -\frac{\wp'(y) \gamma_n'}{(\wp(y) - \gamma_n)^2},\\
d_n(x,y) = \frac{F_1(\gamma_{n-1}) F_1(\gamma_{n}) (\wp(y) - \gamma_{n-2}) (\wp(y) - \gamma_{n+1})}{(\gamma_{n-2} - \gamma_{n-1})(\gamma_{n-1} - \gamma_n)^2 (\gamma_n - \gamma_{n+1})(\wp(y) - \gamma_{n-1})(\wp(y) - \gamma_n)},\\
f_n(x,y) = -\frac{\wp'(y) (\gamma_n - \gamma_{n+1})}{(\wp(y) - \gamma_n) (\wp(y) - \gamma_{n+1})} + g_n(y),\\
g_n(y) = \frac{(-1)^n}{\wp'(y)} \big((n s_1 + s_0) \wp^2(y) + (n k_1 + k_0) \wp(y) + (n p_1 + p_0) \big)
\end{gather*}
are rank two algebro--geometric solutions of (\ref{eq2})--(\ref{eq4}) corresponding to the spectral curve (\ref{eq6}). Here $s_j, k_j, p_j$ are constants, $j=1,2$.
\end{theorem}
Equation (\ref{eq7}) has the following Lax representation
$$
[L_4, \partial_x-V_{n-1}(x)V_n(x)T^{-2}]=0,
$$
where $L_4 = (T+V_n(x)T^{-1})^2 + W_n(x),$
\begin{gather}
\label{eq9}
V_n(x)=\frac{F_1(\gamma_n(x))}{(\gamma_n(x)-\gamma_{n-1}(x))(\gamma_n(x)-\gamma_{n+1}(x))},\\
\label{eq10}
W_n(x)=-c_2-\gamma_n(x)-\gamma_{n+1}(x).
\end{gather}
The operator $L_4$ commutes with an operator $L_6$
and $L_4, L_6$ form a one--point rank two pair of operators with the spectral curve (\ref{eq6}). The equation (\ref{eq7}) can be considered a difference analogue of KN equation
\begin{equation}
\label{eq11}
U_{t}=\frac{48F_1(-\frac{1}{2}(c_2+U))-U_{xx}^2+2U_xU_{xxx}}{8U_x}.
\end{equation}
Equation (\ref{eq11}), as well as (\ref{eq7}), admits the Lax representation related to the rank two (differential) operators corresponding to the elliptic spectral curve.
Moreover, (\ref{eq7}), as well as (\ref{eq11}), appears as an auxiliary equation for the separation variables in the $2+1$ system. For these reasons, we call (\ref{eq7}) the Difference Krichever--Novikov equation (DKN).

Difference chains of type (\ref{eq7}) were studied in many papers (see e.g. \cite{MShYa,LWYa}), but we did not find (\ref{eq7}) described in literature.

In Section 2 we study the rank two algebro--geometric solutions of the system
\begin{gather}
\label{eq12}
\partial_x V_n=V_n(W_{n-1}-W_n+V_{n-1}-V_{n+1}),\\
\label{eq13}
\partial_x W_n=(W_{n}-W_{n-1})V_n+(W_{n+1}-W_n)V_{n+1}.
\end{gather}
This system admits a Lax pair (see (\ref{eq17}) below). This system is reduced to DKN under the reduction (\ref{eq9}), (\ref{eq10}) at $g=1$.

In Section 3 we prove Theorem 1.

\section{DKN Equation}

Let us consider one--point operators of rank two $L_4, L_{4g+2}$ corresponding to the hyperelliptic spectral curve $\Gamma$ given by (\ref{eq5}).
The common eigenfunctions of $L_4$ and $L_{4g+2}$ satisfy the equation
$$
\psi_{n+1}(P)=\chi_1(n,P)\psi_{n-1}(P)+\chi_2(n,P)\psi_n(P),
$$
where $\chi_1(n,P)$ and $\chi_2(n,P)$ are rational functions on $\Gamma$ having
$2g$ simple poles, depending on $n$ (see~\cite{KN2}). The function $\chi_2(n,P)$
has, in addition, a simple pole at $q=\infty$. To find $L_4$ and $L_{4g+2}$ it is sufficient to find $\chi_1$ and $\chi_2.$ Let $\sigma$ be the holomorphic involution on $\Gamma, \ \ \sigma(z,w)=\sigma(z,-w).$
\\
In~\cite{MM1} it was proved that if
$$
\chi_1(n,P)=\chi_1(n,\sigma(P)),\qquad
\chi_2(n,P)=-\chi_2(n,\sigma(P)),
$$
then $L_4$ has the form
\begin{equation}
\label{eq14}
L_4=(T+V_nT^{-1})^2+W_n,
\end{equation}
where
$$
\chi_1=-V_n\frac{Q_{n+1}}{Q_{n}},\qquad
\chi_2=\frac{w}{Q_n}, \qquad
Q_n(z)=z^g+\alpha_{g-1}(n)z^{g-1}+\ldots+\alpha_0(n),
$$
and polynomial $Q$ satisfies the equation
\begin{equation}
\label{eq15}
F_g(z)=Q_{n-1}Q_{n+1}V_n+Q_{n}Q_{n+2}V_{n+1}+Q_nQ_{n+1}(z-V_n-V_{n+1}-W_n).
\end{equation}
From (\ref{eq15}) it follows that $Q$ also satisfies the linear equation
\begin{multline}
\label{eq16}
Q_{n-1}V_n+Q_n(z-V_n-V_{n+1}-W_{n}) - Q_{n+2}(z-V_{n+1}\\ - V_{n+2} - W_{n+1}) - Q_{n+3}V_{n+2}=0.
\end{multline}
At $g=1,$ we have $Q_n=z-\gamma_n,$ and equation (\ref{eq15}) has the solution
$$
V_n=\frac{F_1(\gamma_n)}{(\gamma_n-\gamma_{n-1})(\gamma_n-\gamma_{n+1})}, \qquad
W_n=-c_2-\gamma_n-\gamma_{n+1}.
$$
At $g>1,$ it is a very difficult problem to solve the equation (\ref{eq15}), as it is to find examples of its solutions.

At the end of this section we study difference evolution equations related to the operator $L_4$ (\ref{eq14}).

The system (\ref{eq12}), (\ref{eq13}) has the following Lax representation
%
\begin{equation}
\label{eq17}
[(T+V_n(x)T^{-1})^2+W_n(x),\partial_x-V_{n-1}(x)V_n(x)T^{-2}]=0.
\end{equation}
The system (\ref{eq12}), (\ref{eq13}) is included in the hierarchy of evolution equations of the form
\begin{equation}
\label{eq18}
[(T+V_n(t_k)T^{-1})^2+W_n(t_k),\partial_{t_k}-P_1(n,t_k)T^{-2}-\ldots-P_k(n,t_k)T^{-2 k}]=0.
\end{equation}
These evolution equations define symmetries of (\ref{eq12}), (\ref{eq13}). At $k=2$ we have
\begin{multline}
\label{eq19}
\partial_{t_k} V_n=V_n \big(V_{n-2} V_{n-1} + V_{n-1}V_n - V_n V_{n+1} -V_{n+1} V_{n+2} +
V_{n-1}^2 -
V_{n+1}^2\\ + W_{n-1}^2 - W_n^2 +2 (V_{n-1} + V_n) W_{n-1} - 2 (V_n + V_{n+1}) W_n \big),
\end{multline}
\begin{multline}
\label{eq20}
\partial_{t_k} W_n=V_{n-1} V_n (W_{n-2} - 2 W_{n-1} + W_n) -
V_{n+1} V_{n+2} (W_n - 2 W_{n+1} + W_{n+2})\\ -V_n (W_{n-1} - W_n) (2 V_n + W_{n-1} + W_n)\\ -
V_{n+1} (W_n - W_{n+1}) (2 V_{n+1} + W_n + W_{n+1}).
\end{multline}
In the case of the algebro--geometric operator $L_4$ at $g=1$, i.e. $V_n$ and $W_n$ have the form (\ref{eq9}), (\ref{eq10}) the system (\ref{eq12}), (\ref{eq13}) is reduced to DKN equation and equation (\ref{eq18}) is reduced to an equation from the DKN hierarchy. For example, the system (\ref{eq19}), (\ref{eq20}) is reduced to
\begin{multline*}
\partial_{t_k}\gamma_n=V_n \big(V_{n + 1} (W_{n-1} - 2 W_n + W_{n+1})\\ -
V_{n-1} (W_{n-2} - 2 W_{n-1} + W_n) + (W_{n-1} - W_n) (2 V_n + W_{n-1} + W_n) \big).
\end{multline*}
At $g>1$ there is no explicit reduction of (\ref{eq12}), (\ref{eq13}), since there is no explicit form of $L_4.$ Nevertheless  one can find an evolution equation on polynomial $Q_n$ associated with the algebro--geometric operator $L_4.$ By direct calculation one can check the following lemma.
\begin{lemma}
Equation
\begin{equation}
\label{eq21}
\partial_x Q_n=V_n(Q_{n+1}-Q_{n-1})
\end{equation}
together with (\ref{eq12}), (\ref{eq13}) defines a symmetry of (\ref{eq15}) and (\ref{eq16}).
\end{lemma}
At $g=1$ equation (\ref{eq21}) is equivalent to DKN.

\section{Proof of Theorem 1}

In this section we explain how to obtain a rank two algebro--geometric solution of (\ref{eq2})--(\ref{eq4}) at $g=1.$ A similar method works for KP (see~\cite{Mir,PL}). The main idea is to apply the Darboux type transformation to $L_4.$ If $Q_n$ satisfies (\ref{eq15}) we have the following factorization (see~\cite{MM1})
$$
L_4-z = \big(T+\chi_2(n+1)-\frac{V_{n-1}V_n}{\chi_1(n-1)}T^{-1} \big)\big(T-\chi_2(n)-\chi_1(n)T^{-1} \big).
$$
Let us assume that $\gamma_n=\gamma_n(x)$ and $z=z_0(y).$ After the Darboux transformation we get
$$
\tilde {L}_4 = \big(T-\chi_2(n)-\chi_1(n)T^{-1} \big) \big(T+\chi_2(n+1)-\frac{V_{n-1}V_n}{\chi_1(n-1)}T^{-1} \big)+z_0(y).
$$
Here $V_n=V_n(x)$ has the form (\ref{eq9}),
{\small $$
\chi_1(n) = -V_n(x) \frac{z_0(y)-\gamma_{n+1}(x)}{z_0(y)-\gamma_{n}(x)}, \quad \chi_2(n) = \frac{w(y)}{z_0(y)-\gamma_n(x)},
\quad
w^2(y) = F_1(z_0(y)).
$$}
The operator $\tilde{L}_4$ has the form
\begin{gather*}
\tilde {L}_4 = T^2+A_1(n,x,y)T+A_0(n,x,y)+A_{-1}(n,x,y)T^{-1}+A_{-2}(n,x,y)T^{-2},\\
A_1(n,x,y)=\frac{(\gamma_{n+2}-\gamma_{n})z_0'(y)}{(z_0(y)-\gamma_{n})(z_0(y)-\gamma_{n+2})},\\
A_0(n,x,y) = \frac{V_n (z_0(y)-\gamma_{n+1})^2 + V_{n+1} (z_0(y)-\gamma_{n})^2-F_1(z_0(y))}{(z_0(y)-\gamma_{n}) (z_0(y)-\gamma_{n+1})} + z_0(y),\\
A_{-1}(n,x,y) = \frac{(\gamma_{n-1}-\gamma_{n+1}) V_n z_0'(y)}{(z_0(y)-\gamma_{n})^2},\\
A_{-2}(n,x,y) = \frac{V_{n-1} V_n (z_0(y)-\gamma_{n-2}) (z_0(y)-\gamma_{n+1})}{(z_0(y)-\gamma_{n-1})(z_0(y)-\gamma_{n})}.
\end{gather*}
The operator $\tilde{L}_4$ commutes with $\tilde{L}_6,$ and $\tilde{L}_4, \tilde{L}_6$ are operators of rank two with the same spectral curve (\ref{eq6}).
By direct calculation one can check that if
\begin{gather*}
b_n = -\frac{z_0'(y) \gamma'_n}{(z_0(y) - \gamma_n)^2},\\
d_n = \frac{F_1(\gamma_{n-1}) F_1(\gamma_{n}) (z_0(y) - \gamma_{n-2}) (z_0(y) - \gamma_{n+1})}{(\gamma_{n-2} - \gamma_{n-1})(\gamma_{n-1} - \gamma_n)^2 (\gamma_n - \gamma_{n+1})(z_0(y) - \gamma_{n-1})(z_0(y) - \gamma_n)},
\end{gather*}
and $\gamma_n(x)$ satisfies DKN, then
$$
[\tilde{L}_4,\partial_x-b_n(x,y)T^{-1}-d_n(x,y)T^{-2}] = 0.
$$
By direct calculation one can also check that if $z_0(y)=\wp(y)$ satisfies (\ref{eq8}), $\gamma_{n+4}=\gamma_n$, and $f_n(x,y)$ has the form as in Theorem 1, then
$$
[\tilde{L}_4, \partial_y-T-f_n(x,y)]=0.
$$
Moreover $\tilde{L}_6,$ $\partial_y-T-f_n(x,y),$ $\partial_x-b_n(x,y)T^{-1}-d_n(x,y)T^{-2}$ commute pairwise. Theorem 1 is proved.

Gulnara S.~Mauleshova

Sobolev Institute of Mathematics,
 	
4 Acad. Koptyug avenue, 630090, Novosibirsk, Russia

${\ }$

Novosibirsk State University,

1, Pirogova str., 630090, Novosibirsk, Russia

e-mail: mauleshova\_gs@mail.ru

${\ }$

${\ }$ ${\ }$

Andrey E.~Mironov

Sobolev Institute of Mathematics,

4 Acad. Koptyug avenue, 630090, Novosibirsk, Russia

${\ }$

Novosibirsk State University,

1, Pirogova str., 630090, Novosibirsk, Russia

e-mail: mironov@math.nsc.ru

\end{document}